\documentclass{amsart}
\usepackage{amsmath,amsthm, amscd, amssymb, amsfonts}

\numberwithin{equation}{section}
\begingroup

\newtheorem{theor}[equation]{Theorem}
\newtheorem{lemma}[equation]{Lemma}
\theoremstyle{definition}
\newtheorem{defin}[equation]{Definition}
\newtheorem{remar}[equation]{Remark}
\endgroup

\newcommand{\acad}{\makebox[12pt]{$\circ$\hspace{-2.5pt}\raisebox{1.6pt}{$\scriptscriptstyle{\succ}$}}}
\newcommand{\acai}{\makebox[12pt]{\raisebox{1.4pt}{$\scriptscriptstyle{\prec}$}\hspace{-2.5pt}$\circ$}}
\newcommand{\aaad}{\makebox[12pt]{$\ast$\hspace{-2pt}\raisebox{1.4pt}{$\scriptscriptstyle{\succ}$}}}

\newcommand{\al}{\alpha}
\newcommand{\be}{\beta}

\newcommand{\CC}{\mathbb{C}}
\newcommand{\Ss}{\mathbb{S}}
\newcommand{\Aut}{\operatorname{Aut}}
\newcommand{\Der}{\operatorname{Der}}
\newcommand{\id}{\operatorname{id}}
\newcommand{\adj}{\operatorname{adj}}

\newcommand{\ve}{\varepsilon}
\newcommand{\Oo}{\mathcal{O}}
\newcommand{\lieg}{\mathfrak{g}}
\newcommand{\lieh}{\mathfrak{h}}
\newcommand{\liegt}{\widetilde{\mathfrak{g}}}
\newcommand{\liea}{\mathfrak{a}}
\newcommand{\vect}{\operatorname{Vect}}
\newcommand{\spn}{\operatorname{span}}
\newcommand{\im}{\operatorname{Im}}
\newcommand{\rhoga}{\rho_{\lieg\liea}}
\newcommand{\rhogax}{\rho_{\lieg\liea X}}
\newcommand{\pmu}{\pi^{-1}}
\newcommand{\mdpd}[4]{\left(\begin{smallmatrix}#1&#2\\#3&#4\end{smallmatrix}\right)}
\newcommand{\mtpt}[9]{\left(\begin{smallmatrix}#1&#2&#3\\#4&#5&#6\\#7&#8&#9\\\end{smallmatrix}\right)}

\begin{document}
\title{Quantization of non-unitary geometric classical $r$-matrices}
\author{P. Etingof \& M. Gra\~na}
\begin{abstract}
In this paper we explicitly attach to a \emph{geometric} 
classical r-matrix $r$ (not necessarily unitary),
a geometric (i.e., set-theoretical) quantum R-matrix $R$, which is a quantization of 
$r$. To accomplish this, we use 
the language of bijective cocycle $7$-tuples, developed by A. Soloviev
in the study of set-theoretical quantum R-matrices.
Namely, we define a classical version of bijective cocycle $7$-tuples, and show
that there is a bijection between them and geometric classical r-matrices.
Then we show how any classical bijective cocycle $7$-tuple can be quantized,
and finally use Soloviev's construction, which turns a 
(quantum) bijective cocycle $7$-tuple into a geometric quantum R-matrix. 
\end{abstract}
\address{Pavel Etingof:\newline\indent
        MIT Math. Dept. Of. 2-176 \newline\indent
        77 Mass. Ave. \newline\indent
        02139, Cambridge, MA, USA.}
\email{etingof@math.mit.edu}
\address{Mat\'\i as Gra\~na:\newline\indent
        MIT Math. Dept. Of. 2-155\newline\indent
        77 Mass. Ave.\newline\indent
        02139, Cambridge, MA, USA.\newline\indent
		\indent Permanent: \newline\indent
		\indent Depto. Matem\'atica - FCEyN - UBA.\newline\indent
		\indent Pab. I - Ciudad Universitaria.\newline\indent
		\indent (1428) Ciudad de Buenos Aires - Argentina}
\email{matiasg@math.mit.edu}
\thanks{The work of P.E. was supported by the NSF grant 
DMS-9988796. The work of M.G. was supported by Conicet.}
\maketitle

\section{Introduction}
Let $X$ be a smooth affine algebraic variety over $\CC$. 
A formal diffeomorphism of $X$ is an automorphism
of the ring $\Oo[X][[\hbar]]$ which is the identity modulo $\hbar$. That is, it looks
like $1+\hbar r+O(\hbar^2)$. Note that for such a map to be a ring homomorphism, $r$
must be a derivation, i.e., $r\in\vect(X)$.

A geometric solution of the Quantum Yang--Baxter Equation (QYBE) is a formal
diffeomorphism $R$ of $X\times X$ such that $R^{12}R^{13}R^{23}=R^{23}R^{13}R^{12}$
as formal diffeomorphisms of $X\times X\times X$. It is straightforward to
check that if $R=1+\hbar r+O(\hbar^2)$ is a geometric solution of the QYBE,
then $r$ is a geometric solution of the Classical Yang--Baxter Equation (CYBE),
i.e. $[r^{12},r^{13}]+[r^{12},r^{23}]+[r^{13},r^{23}]=0$.
The well known \emph{quantization problem} is: given
a geometric solution $r$ of the CYBE, construct a geometric solution 
$R$ of the QYBE which restricts to $r$ in degree $1$, i.e., $R=1+\hbar r+O(\hbar^2)$.
Such an $R$ is called a \emph{geometric quantization} of $r$.

In \cite{ek} it is proved that any solution of CYBE can be quantized.
However, the proof does not give a simple explicit formula for $R$,
and furthermore it is not clear if the map $R$ is geometric when $r$ is.
On the other hand, in \cite{es}, the problem of geometric quantization is
solved for unitary r-matrices (i.e., satisfying the unitarity condition
$r^{21}=-r$): it is shown that in this case the geometric quantization $R$
exists and satisfies the quantum unitarity condition $R^{21}R=1$. For this,
using the approach of \cite{ess}, it is shown that both classical and
quantum geometric r-matrices are parametrized by some group-theoretical
data, at the level of which quantization basically reduces to the usual
exponential map.

In this paper we generalize the result of \cite{es} to the non-unitary
case.  Namely, we prove the following result.

\begin{theor}\label{tth}
Any geometric solution to the CYBE admits a geometric quantization.
\end{theor}

As in the unitary case, this is accomplished using the group-theoretical approach, 
developed in \cite{s}.  In this sense, 
this paper completes
the square
$$\begin{CD}
\mbox{\cite{ess}} @>\mbox{forget unitarity}>> \mbox{\cite{s}} \\
@V\mbox{\begin{minipage}{1.5cm}
	\begin{flushright}geometric \\[-6pt] version\end{flushright}
\end{minipage}}VV
@VV\mbox{\begin{minipage}{1.5cm}
	\begin{flushleft}geometric \\[-6pt] version\end{flushleft}
\end{minipage}}V \\
\mbox{\cite{es}} @>\mbox{forget unitarity}>>
	\mbox{\begin{minipage}{1cm}\begin{center}this\\[-3pt]paper\end{center}\end{minipage}}
\end{CD}
$$

\section{Cocycle $7$-tuples}\label{ss2}

Recall the following definition from \cite{s}.
\begin{defin} 
A \emph{bijective cocycle $7$-tuple} (BCST) is a $7$-tuple
$$(G,A,X,\rho_{GA},\rho_{GAX},\pi,\Psi),$$
where $G$ and $A$ are groups, $X$ a set, $\rho_{GA}$ is an action of $G$ on $A$,
$\rho_{GAX}$ is an action of $A\rtimes G$ on $X$, $\pi:G\to A$
is a bijective $1$-cocycle (i.e., $\pi(gh)=\pi(g)(g\pi(h))$)
and $\Psi:X\to A$ is an $A\rtimes G$-equivariant map (with $A$ acting on itself
by conjugation) whose image generates $A$.

A BCST gives two actions of $G$ on $X$, given by composing $\rho_{GAX}$ with
the inclusion $G\to A\rtimes G$ and with the map $g\mapsto (\pi(g), g)$.
Thus, it defines a map $G\to\Ss_X\times\Ss_X$ (where $\Ss_X$ is the group
of bijections $X\to X$).
The BCST is said to be \emph{faithful} if this map is injective.
\end{defin}

We state a result from \cite{s}, slightly modified to fit our definition
of $1$-cocycles:

\begin{lemma}\label{lm:sfst}
Let $(G,A,X,\rho_{GA},\rho_{GAX},\pi,\Psi)$, be a BCST.
Let $R:X\times X\to X\times X$ be defined by
\begin{equation}\label{eq:sos}
R(x,y)=(x\check *y, x\check\circ y),
\end{equation}
where $x\check *y=\rho_{GAX}(\pi^{-1}(y^{-1}))(x)$ and
$x\check\circ y=\rho_{GAX}(\pi^{-1}((x\check *y)^{-1}))^{-1}(\rho_{GAX}(x\check *y)(y))$
(we denote here for $z\in X$, $z$ instead of $\Psi(z)\in A$ for not
overcharging the notation).
Then $R$ is a set theoretical solution of the braid equation, i.e.,
$R^{12}R^{13}R^{23}=R^{23}R^{13}R^{12}$.
\end{lemma}

In fact, the solution \eqref{eq:sos} is also \emph{non-degenerate} in the sense of
\cite{s}. The main result of \cite{s} is that 
formula \eqref{eq:sos}
defines a bijection (or, more precisely, a categorical equivalence)
between faithful bijective cocycle 7-tuples and nondegenerate set-theoretical solutions 
of QYBE.

Let us now generalize this theory to the quasiclassical situation. 
Our ground field will always be $\CC$.

\begin{defin}
A \emph{classical bijective cocycle $7$-tuple} (CBCST) is a $7$-tuple
$$(\lieg,\liea,X,\rhoga,\rhogax,\pi,\Psi),$$ where $\lieg$ and $\liea$ are
Lie algebras, $X$ a smooth affine algebraic variety, $\rhoga$ is an action
of $\lieg$ on $\liea$, $\rhogax$ is an action of $\liea\rtimes\lieg$ on
$X$, $\pi:\lieg\to\liea$ is a bijective (non-commutative) $1$-cocycle and
$\Psi:X\to\liea$ is an $\liea\rtimes\lieg$-equivariant map (with $\liea$
acting on itself by commutator) whose image generates $\liea$.

A CBCST gives two actions of $\lieg$ on $X$, given by composing $\rhogax$ with
the inclusion $\lieg\to\liea\rtimes\lieg$ and with the map $g\mapsto (\pi(g), g)$.
Thus, it defines a map $\lieg\to\vect(X)\times\vect(X)$.
The CBCST is said to be \emph{faithful} if this map is injective.
\end{defin}

The first part of the proof of Theorem \ref{tth} is to relate classical bijective
cocycle $7$-tuples and geometric solutions of the CYBE.
One can consider the categories of CBCST's and of geometric solutions of the CYBE,
with the obvious notion of morphisms in both.

\begin{theor}\label{tle}
There is a an equivalence of categories between faithful CBCST's and geometric solutions
of the CYBE.
\end{theor}

\begin{proof}
We will construct 
mutually (quasi)inverse functors in both directions.
More precisely, we will do this only for objects, as the extension to 
morphisms is straightforward (and not used).
In the proof, we will refer to a sequence of lemmas, which are
stated and proved in \S3. 

Let $r$ be a geometric solution of the CYBE. Then we have
\begin{align}
\nonumber		&r \in \vect(X)\otimes\Oo(X)\;\oplus\;\Oo(X)\otimes\vect(X)\;, \\
\label{eq:trm}	&r=\sum_ia^1_i\otimes a^0_i+\sum_jb^0_j\otimes b^1_j\;, \\
\nonumber		&[r^{12},r^{13}]+[r^{12},r^{23}]+[r^{13},r^{23}]=0\;.
\end{align}
We will omit the summation sign in expressions of the type
$\sum_ia^1_i\otimes a^0_i$ or $\sum_jb^0_j\otimes b^1_j$.
By looking at the degree of components, we can split the CYBE into
three equations as follows:
\begin{align}
\label{eq:cybe1}0=\;&[a^1_i,a^1_k]\otimes a^0_i\otimes a^0_k
	\;-\;a^1_i\otimes a^1_k\cdot a^0_i\otimes a^0_k
	\;-\;a^1_i\otimes b^0_l\otimes b^1_l\cdot a^0_i\;, \\
\label{eq:cybe2}0=\;&-a^1_k\cdot b^0_j\otimes b^1_j\otimes a^0_k
	\;+\;b^0_j\otimes[b^1_j,a^1_k]\otimes a^0_k
	\;+\;b^0_j\otimes a^1_k\otimes b^1_j\cdot a^0_k\;, \\
\label{eq:cybe3}0=\;&a^1_i\cdot b^0_l\otimes a^0_i\otimes b^1_l
	\;+\;b^0_j\otimes b^1_j\cdot b^0_l\otimes b^1_l
	\;+\;b^0_j\otimes b^0_l\otimes [b^1_j,b^1_l]\;.
\end{align}

We define $\lieg_1=\spn\{a^1_i\}_i$, $\lieg_2=\spn\{b^1_j\}_j$
(we assume that the expression \eqref{eq:trm} has the minimal
possible number of summands).
It is easy to see from \eqref{eq:cybe1} and \eqref{eq:cybe3}
that both $\lieg_i$'s are Lie subalgebras of $\vect(X)$.
Also, $[\lieg_1,\lieg_2]\subseteq\lieg_1+\lieg_2$ by \eqref{eq:cybe2}.
For $x\in\Oo(X)^*$, let
$$p(x)=(-x(a^0_i)a^1_i, x(b^0_j) b^1_j) \in \lieg_1\oplus\lieg_2,$$
and call $\lieg=\im(p)$. It is proved in \eqref{eq:fryb} below that $\lieg$ is
a Lie subalgebra of $\lieg_1\oplus\lieg_2$.

We define $V_1,\;V_2\subset\Oo(X)$ as $V_1=\spn\{a^0_i\}_i$, $V_2=\spn\{b^0_j\}_j$,
and we put $\liea=(V_1+V_2)^*$. We can restrict $p$ to $\liea$ and we get an
isomorphism $p|_{\liea}:\liea\to\lieg$. Call $\pi=p|_{\liea}^{-1}$.

Consider the action of $\vect(X)$ on $\Oo(X)$.
From \eqref{eq:cybe1}, \eqref{eq:cybe2}, \eqref{eq:cybe3} one can see that
\begin{align*}
&\lieg_1\cdot V_1\subset V_1+V_2,\quad \lieg_1\cdot V_2\subset V_2\;, \\
&\lieg_2\cdot V_1\subset V_1,\quad \lieg_2\cdot V_2\subset V_1+V_2\;.
\end{align*}

We denote by $\aaad$ and $\acai$ the left and right actions of $\lieg$ on $\liea$ via
projections to the first and second coordinates. We also use this notation
for the maps $\liea\times\liea\to\liea$ 
obtained by composing these actions with $p$ (we warn that these maps are
not actions since $p$ is not a Lie algebra homomorphism). Specifically,
\begin{align}
\label{eq:star}
(x\aaad y)(f)&=(p(x)\aaad y)(f)=-x(a^0_i)(a^1_i\cdot y)(f)=x(a^0_i)y(a^1_i\cdot f)\;, \\
\label{eq:circ}
(x\acai y)(f)&=(x\acai p(y))(f)=y(b^0_j)(x\cdot b^1_j)(f)=y(b^0_j)x(b^1_j\cdot f)\;.
\end{align}
Using these actions and the bijection $\pi$
we equip $\liea$ with a Lie algebra structure:
\begin{equation}\label{eq:dba}
[x,y]=- x\aaad y+ y\aaad x+\pi([p(x),p(y)])\;.
\end{equation}
We have to prove that this is indeed a Lie algebra structure, which we do
in Lemma \ref{lm:iila} below. Then $\pi$ is automatically a $1$-cocycle.

We prove in Lemma \ref{lm:iiad} that the action $\aaad $ of $\lieg$ on $\liea$ is by
derivations (i.e., $\aaad :\lieg\to\Der(\liea)$). This allows us to take the semidirect
product $\liea\rtimes\lieg$, whose structure, we recall, is
$$[(a,g),(b,h)]=(g\aaad b-h\aaad a+[a,b],[g,h])\;.$$

Notice from \eqref{eq:sryb} below that the Lie algebra structure in $\liea$ is
\begin{equation}\label{eq:lsa}
[x,y]=-y(a^0_i)(a^1_i\cdot x)-y(b^0_j)(b^1_j\cdot x)
	=x(a^0_i)(a^1_i\cdot y)+x(b^0_j)(b^1_j\cdot y)\;.
\end{equation}

Consider the action of $\lieg$ on $X$ given by the projection to the
first coordinate, i.e.
\begin{equation}\label{eq:rgx}
\rho(p(x))=-x(a^0_i)a^1_i\;.
\end{equation}
Let $d:\lieg\to\lieg_1+\lieg_2$, $d(g_1,g_2)=g_1-g_2$.
Define $\rho_{\liea X}:\liea\to\vect(X)$ as follows:

\begin{equation}\label{eq:rax}
\begin{CD}
\rho=\Big(\liea @>-p>> \lieg @>d>> \lieg_1+\lieg_2 \subset \vect(X)\Big),\quad
\text{i.e. } \rho(x)=x(a^0_i)a^1_i+x(b^0_j)b^1_j\;.
\end{CD}
\end{equation}
We prove in Lemma \ref{lm:iiam} that this is a map of Lie algebras.
We can lift this action to $\liea\rtimes\lieg$: define
$$\rho_{\lieg\liea X}:\liea\rtimes\lieg\to\vect X,\quad
	\rho(y,p(x))=y(a^0_i)a^1_i+y(b^0_j)b^1_j-x(a^0_k)a^1_k\;.$$
From the structure of $\liea\rtimes\lieg$ we have $[(0,p(x)),(y,0)]=(x\aaad y,0)$.
Let us check that this is indeed an action. We have
\begin{align*}
\rho(x\aaad y) &= (x\aaad y)(a^0_i)a^1_i+(x\aaad y)(b^0_j)b^1_j
	=x(a^0_k)y(a^1_k\cdot a^0_i)a^1_i+x(a^0_k)y(a^1_k\cdot b^0_j)b^1_j\;, \\
[\rho(p(x)),\rho_{\liea X}(y)] &= [-x(a^0_k)a^1_k,y(a^0_i)a^1_i+y(b^0_j)b^1_j]
	=x(a^0_k)y(a^0_i)[a^1_i,a^1_k]-x(a^0_k)y(b^0_j)[a^1_k,b^1_j] \\
	&=y(a^1_k\cdot a^0_i)x(a^0_k)a^1_i+y(b^0_l)x(b^1_l\cdot a^0_i)a^1_i
		+y(a^1_k\cdot b^0_j)x(a^0_k)b^1_j-y(b^0_j)x(b^1_j\cdot a^0_k)a^1_k \\
	&=\rho(x\aaad y)\;,
\end{align*}
as desired. 

Define $\Psi_{\liea}:X\to\liea$ by restriction: $\Psi_{\liea}(x)(f)=f(x)$ for $f\in V$.
Comparing \eqref{eq:lsa} with \eqref{eq:rax} and \eqref{eq:star} with \eqref{eq:rgx},
it is clear that $\Psi_{\liea}$ is $\liea\rtimes\lieg$-invariant, i.e.,
$$\Psi_{\liea *}|_x\rho_{\lieg\liea X}(a,g)=-(a,g)\cdot\Psi_{\liea}(x)
	=-[a,\Psi_{\liea}(x)]-g\aaad \Psi_{\liea}(x).$$
We have thus constructed a CBCST $(\lieg, \liea, X, \rhoga, \rhogax, \pi, \Psi_{\liea})$.

\bigskip
Now we begin with a $7$-tuple $(\lieg,\liea,X,\rhoga,\rhogax,\pi,\Psi_{\liea})$
and aim to construct $r$. First, consider the maps
\begin{align*}
&\begin{CD}
\alpha =
X @>\Psi_\liea>> \liea @>-\pmu>> \lieg @>i>> \liea\rtimes\lieg @>\rho>> \vect X\;,
\end{CD} \\
&\begin{CD}
\beta =
X @>\Psi_\liea>> \liea @>\pmu\times\id >> \liea\rtimes\lieg @>\rho>> \vect X\;.
\end{CD} \\
\end{align*}
We call $\lieg_1=\im(\alpha)$, $\lieg_2=\im(\beta)$. Now, $\al$ gives
by composition a map $\lieg_1^*\to\Oo(X)$, which in turn is an element of
$\lieg_1\otimes\Oo(X)$. Analogously, $\be$ gives a map $\lieg_2^*\to\Oo(X)$, and in
turn an element of $\Oo(X)\otimes\lieg_2$. Call these elements $r^1,r^2$ respectively.
We view them as elements of $\vect(X)\otimes\Oo(X)$ and $\Oo(X)\otimes\vect(X)$. Call
\begin{equation}\label{eq:raas}
r=r^1+r^2.
\end{equation}
We prove in Lemma \ref{lm:iias} that $r$ is a geometric solution of the CYBE.
It is easy to see that both constructions
$$\text{Geometric solutions of CYBE} \leftrightsquigarrow
	\text{Classical bijective cocycle $7$-tuples}$$
are inverse to each other. The condition on $\Psi(X)$ to generate $\liea$ guarantees
that after applying CBCST $\rightsquigarrow$ Geom. sol. of CYBE $\rightsquigarrow$ CBCST
one gets an algebra $\liea$ isomorphic to the original one. Analogously, the faithfulness
condition guarantees that one recovers an algebra $\lieg$ isomorphic to the original one.
\end{proof}

\section{Auxiliary lemmas}

\begin{lemma}\label{lm:iila}
The definition in \eqref{eq:dba} equips $\liea$ with the
structure of a Lie algebra.
\end{lemma}
\begin{proof}
The bilinearity and antisymmetry of the bracket are clear. Let us prove
that it satisfies the Jacobi identity. We have
\begin{align*}
p(x\aaad y) &= (-(x\aaad y)(a^0_i)a^1_i,(x\aaad y)(b^0_j)b^1_j) \\
	&= (-x(a^0_k)y(a^1_k\cdot a^0_i)a^1_i,x(a^0_k)y(a^1_k\cdot b^0_j)b^1_j)\;, \\
p(x\acai y) &= (-(x\acai y)(a^0_i)a^1_i,(x\acai y)(b^0_j)b^1_j) \\
	&= (-y(b^0_l)x(b^1_l\cdot a^0_i)a^1_i,y(b^0_l)x(b^1_l\cdot b^0_j)b^1_j)\;, \\
[p(x),p(y)] &= [(-x(a^0_i)a^1_i,x(b^0_j) b^1_j),(-y(a^0_k)a^1_k,y(b^0_l) b^1_l)] \\
	&= ([x(a^0_i)a^1_i,y(a^0_k)a^1_k], [x(b^0_j) b^1_j,y(b^0_l) b^1_l]) \\
	&= (x(a^0_i)y(a^0_k)[a^1_i,a^1_k],x(b^0_j)y(b^0_l)[b^1_j,b^1_l]) \\
	&= (-x(a^0_k)y(a^0_i)[a^1_i,a^1_k],-x(b^0_l)y(b^0_j)[b^1_j,b^1_l])\;, \\
\end{align*}
whence, using \eqref{eq:cybe1} and \eqref{eq:cybe3},
\begin{align}
\label{eq:fryb}
[p(x),p(y)]&=p(x\aaad y)+p(x\acai y),\quad\text{i.e.} \\
\label{eq:sryb}
[x,y]&=x\acai y+y\aaad x\;.
\end{align}
Notice that \eqref{eq:fryb} proves that $\lieg$ is a subalgebra of $\lieg_1\oplus\lieg_2$.
Similarly we compute for $f\in\Oo(X)$,
\begin{align*}
((z\aaad y)\acai x)(f) &= x(b^0_j)(z\aaad y)(b^1_j\cdot f)
	=x(b^0_j)y(a^1_k b^1_j\cdot f)z(a^0_k)\;, \\
(z\aaad (y\acai x))(f) &= z(a^0_i)(y\acai x)(a^1_i\cdot f)
	=x(b^0_l)y(b^1_l a^1_i\cdot f)z(a^0_i)\;, \\
((z\acai x)\aaad y)(f) &= (z\acai x)(a^0_i)y(a^1_i\cdot f)
	=x(b^0_l)y(a^1_i\cdot f)z(b^1_l\cdot a^0_i)\;, \\
(y\acai (z\aaad x))(f) &= (z\aaad x)(b^0_j)y(b^1_j\cdot f)
	=x(a^1_k\cdot b^0_j)y(b^1_j\cdot f)z(a^0_k)\;.
\end{align*}
Using \eqref{eq:cybe2} we get
\begin{align}
\nonumber &((z\aaad y)\acai x)-(z\aaad (y\acai x))-((z\acai x)\aaad y)+(y\acai (z\aaad x))=0\;, \\
\intertext{and using now \eqref{eq:fryb}, we have}
\label{eq:et}
&\pi[p(z\aaad y),p(x)]-((z\aaad y)\aaad x)-(z\aaad \pi[p(y),p(x)])+(z\aaad (y\aaad x)) \\
\nonumber&\qquad-(\pi[p(z),p(x)]\aaad y)+((z\aaad x)\aaad y)+\pi[p(y),p(z\aaad x)]-(y\aaad (z\aaad x))=0.
\end{align}
Now, we compute
\begin{align*}
[x,[y,z]]
=\;& -x\aaad [y,z]+[y,z]\aaad x+\pi[p(x),p([y,z])] \\
	=\;& x\aaad (y\aaad z)-x\aaad (z\aaad y)-x\aaad (\pi[p(y),p(z)])-(y\aaad z)\aaad x \\
	&+(z\aaad y)\aaad x+\pi[p(y),p(z)]\aaad x-\pi[p(x),p(y\aaad z)]\\
	&+\pi[p(x),p(z\aaad y)]+\pi[p(x),[p(y),p(z)]]\;, \\
[y,[z,x]]
	=\;& y\aaad (z\aaad x)-y\aaad (x\aaad z)- y\aaad (\pi[p(z),p(x)]) \\
	&-(z\aaad x)\aaad y+(x\aaad z)\aaad y+\pi[p(z),p(x)]\aaad y \\
	&-\pi[p(y),p(z\aaad x)]+\pi[p(y),p(x\aaad z)]+\pi[p(y),[p(z),p(x)]]\;, \\
[z,[x,y]]
	=\;& z\aaad (x\aaad y)-z\aaad (y\aaad x)- z\aaad (\pi[p(x),p(y)]) \\
	&-(x\aaad y)\aaad z+(y\aaad x)\aaad z+\pi[p(x),p(y)]\aaad z \\
	&-\pi[p(z),p(x\aaad y)]+\pi[p(z),p(y\aaad x)]+\pi[p(z),[p(x),p(y)]]\;.
\end{align*}
Applying \eqref{eq:et} three times and Jacobi identity in $\lieg$ we get
the Jacobi identity in $\liea$.
\end{proof}

\begin{lemma}\label{lm:iiad}
The algebra $\lieg$ acts on $\liea$ by derivations
with the action $\aaad $ defined in \eqref{eq:star}.
\end{lemma}
\begin{proof}
This is straightforward: we compute
\begin{align*}
p(x\aaad [y,z]) &= p(x\aaad (z\aaad y)+x\aaad (y\acai z)) \\
	&=(-x(a^0_i)(z\aaad y+y\acai z)(a^1_i\cdot a^0_k)a^1_k,
		x(a^0_i)(z\aaad y+y\acai z)(a^1_i\cdot b^0_l)b^1_l) \\
	&=(-x(a^0_i)z(a^0_m)y(a^1_ma^1_i\cdot a^0_k)a^1_k
			-x(a^0_i)z(b^0_n)y(b^1_na^1_i\cdot a^0_k)a^1_k, \\
	&\qquad x(a^0_i)z(a^0_m)y(a^1_ma^1_i\cdot b^0_l)b^1_l
			+x(a^0_i)z(b^0_n)y(b^1_na^1_i\cdot b^0_l)b^1_l)\;, \\
p(-[x\aaad y,z])
	&=(((x\aaad y)\acai z)(a^0_k)a^1_k+(z\aaad (x\aaad y))(a^0_k)a^1_k, \\
		&\qquad -((x\aaad y)\acai z)(b^0_l)b^1_l-(z\aaad (x\aaad y))(b^0_l)b^1_l) \\
	&=(z(b^0_j)(x\aaad y)(b^1_j\cdot a^0_k)a^1_k+z(a^0_i)(x\aaad y)(a^1_i\cdot a^0_k)a^1_k, \\
		&\qquad -z(b^0_j)(x\aaad y)(b^1_j\cdot b^0_l)b^1_l
			-z(a^0_i)(x\aaad y)(a^1_i\cdot b^0_l)b^1_l) \\
	&=(z(b^0_j)x(a^0_m)y(a^1_mb^1_ja^0_k)a^1_k
			+z(a^0_i)x(a^0_m)y(a^1_ma^1_ia^0_k)a^1_k, \\
		&\qquad -z(b^0_j)x(a^0_m)y(a^1_mb^1_jb^0_l)b^1_l
			-z(a^0_i)x(a^0_m)y(a^1_ma^1_ib^0_l)b^1_l)\;, \\
p(-[y,x\aaad z])
	&=((y\acai(x\aaad z))(a^0_k)a^1_k+((x\aaad z)\aaad y)(a^0_k)a^1_k, \\
		&\qquad -(y\acai(x\aaad z))(b^0_l)b^1_l-((x\aaad z)\aaad y)(b^0_l)b^1_l) \\
	&=((x\aaad z)(b^0_j)y(b^1_j\cdot a^0_k)a^1_k+(x\aaad z)(a^0_i)y(a^1_i\cdot a^0_k)a^1_k, \\
		&\qquad -(x\aaad z)(b^0_j)y(b^1_j\cdot b^0_l)b^1_l
			-(x\aaad z)(a^0_i)y(a^1_i\cdot b^0_l)b^1_l) \\
	&=(x(a^0_m)z(a^1_m\cdot b^0_j)y(b^1_j\cdot a^0_k)a^1_k
			+x(a^0_m)z(a^1_m\cdot a^0_i)y(a^1_i\cdot a^0_k)a^1_k, \\
		&\qquad -x(a^0_m)z(a^1_m\cdot b^0_j)y(b^1_j\cdot b^0_l)b^1_l
			-x(a^0_m)z(a^1_m\cdot a^0_i)y(a^1_i\cdot b^0_l)b^1_l)\;, \\
\end{align*}
and now \eqref{eq:cybe1} and \eqref{eq:cybe2} apply to see that
$p(x\aaad [y,z]-[x\aaad y,z]-[y,x\aaad z])=0$.
\end{proof}

\begin{lemma}\label{lm:iiam}
The map $\rho_{\liea X}$ in \eqref{eq:rax} is a Lie algebra homomorphism.
\end{lemma}
\begin{proof}
We compute
\begin{align*}
\rho([x,y]) &= \rho(x\acai y+y\aaad x) \qquad\text{(by \eqref{eq:sryb})} \\
	&=(x\acai y+y\aaad x)(a^0_i)a^1_i+(x\acai y+y\aaad x)(b^0_j)b^1_j \\
	&=y(b^0_l)x(b^1_l\cdot a^0_i)a^1_i+y(a^0_k)x(a^1_k\cdot a^0_i)a^1_i
		+y(b^0_l)x(b^1_l\cdot b^0_j)b^1_j+y(a^0_k)x(a^1_k\cdot b^0_j)b^1_j\;, \\
[\rho(x),\rho(y)] &= x(a^0_i)y(a^0_k)[a^1_i,a^1_k]+x(a^0_i)y(b^0_l)[a^1_i,b^1_l]
		+x(b^0_j)y(a^0_k)[b^1_j,a^1_k]+x(b^0_j)y(b^0_l)[b^1_j,b^1_l] \\
	&=x(a^1_k\cdot a^0_i)y(a^0_k)a^1_i+x(b^0_l)y(b^1_l\cdot a^0_i)a^1_i
			-y(a^1_k\cdot b^0_j)x(a^0_k)b^1_j+y(b^0_j)x(b^1_j\cdot a^0_k)a^1_k \\
		&\qquad +x(a^1_k\cdot b^0_j)y(a^0_k)b^1_j -x(b^0_j)y(b^1_j\cdot a^0_k)a^1_k
			-x(a^1_i\cdot b^0_l)y(a^0_i)b^1_l-x(b^0_j)y(b^1_j\cdot b^0_l)b^1_l \\
	&=x(a^1_k\cdot a^0_i)y(a^0_k)a^1_i-y(a^1_k\cdot b^0_j)x(a^0_k)b^1_j
		+y(b^0_j)x(b^1_j\cdot a^0_k)a^1_k-x(b^0_j)y(b^1_j\cdot b^0_l)b^1_l\;,
\end{align*}
and, clearly,
\begin{multline*}
y(b^0_l)x(b^1_l\cdot b^0_j)b^1_j+y(a^0_k)x(a^1_k\cdot b^0_j)b^1_j
	+y(a^1_k\cdot b^0_j)x(a^0_k)b^1_j+x(b^0_j)y(b^1_j\cdot b^0_l)b^1_l \\
	=([x,y]+[y,x])(b^0_j)b^1_j=0\;.
\end{multline*}
\end{proof}

\begin{lemma}\label{lm:iias}
The map $r$ defined in \eqref{eq:raas} is a geometric solution of the CYBE.
\end{lemma}
\begin{proof}
We must prove that $r$ satisfies \eqref{eq:cybe1}, \eqref{eq:cybe2} and \eqref{eq:cybe3}.
To prove \eqref{eq:cybe1} we evaluate the second and third tensorand in each term at
points $b,c\in X$. Let $\liegt=\lieg_1+\lieg_2$ and let $\{x_i\}$, $\{x^i\}$ be dual bases
of $\liegt$ and $\liegt^*$. In order to make formulas more readable, we call
$\{y_j\}$, $\{y^j\}$ another copy of the dual bases.

We identify vector spaces with their tangent spaces. We have:

\begin{align*}
[a^1_i,a^1_k]\otimes a^0_i\otimes a^0_k \mapsto&
	[x_i,y_j] x^i(\al(b)) y^j(\al(c))=[\al(b),\al(c)]\;, \\
-a^1_i\otimes a^1_k\cdot a^0_i\otimes a^0_k \mapsto&
	 -x_i (y_j\cdot (x^i\al))(b)(y^j\al)(c)
	=-x_i (\al(c)\cdot (x^i\al))(b) \\
	&=-x_i \al_{*b}(\al(c))(x^i)
	=-\al_{*b}(\al(c))\;, \\
-a^1_i\otimes b^0_l\otimes b^1_l\cdot a^0_i \mapsto&
	 -x_i (y^j\be)(b)(y_j\cdot(x^i\al))(c)
	=\cdots
	=-\al_{*c}(\be(b))\;.
\end{align*}
Now, $\Psi_\liea$ is $\liea\rtimes\lieg$-invariant, which means that
\begin{align*}
-\al_{*b}(\al(c)) &
	=-(\rho i\pmu)_*(\Psi_\liea)_*|_b(\rho(i\pmu\Psi_\liea(c)))
		=(\rho i\pmu)_*((i\pmu\Psi_\liea(c)) \cdot  (\Psi_\liea)(b)) \\
	&=\rho(\pmu((\pmu\Psi_\liea(c)) \cdot  (\Psi_\liea(b))),0)\;. \\
\end{align*}
Analogously, setting for brevity $\tilde b=\Psi_\liea b$, $\tilde c=\Psi_\liea c$,
we have
\begin{align*}
-\al_{*c}(\be(b))
	&=-\rho\big(\pmu((\pmu\tilde b)\cdot\tilde c+[\tilde b,\tilde c]),0\big)\;, \\
[\al(b),\al(c)]
	&=\rho([\pmu\tilde b,\pmu\tilde c],0)\;,
\end{align*}
and thus
\begin{align*}
\text{\eqref{eq:cybe1}}
	&=\rho\Big(\pmu((\pmu\tilde c) \cdot \tilde b)
		-\pmu((\pmu\tilde b)\cdot(\tilde c)
		-[\tilde b,\tilde c])
		+[\pmu\tilde b,\pmu\tilde c],0\Big) \\
	&=\rho\Big(-\pmu\big(-(\pmu\tilde c) \cdot \tilde b
		+(\pmu\tilde b)\cdot(\tilde c)
		+[\tilde b,\tilde c]\big)
		+[\pmu\tilde b,\pmu\tilde c],0\Big) \\
	&=0\qquad\text{by \eqref{eq:dba}}\;.
\end{align*}

Analogously, evaluating \eqref{eq:cybe2} at $a,c$ and \eqref{eq:cybe3} at $a,b$, we get
\begin{align*}
-a^1_k\cdot b^0_j\otimes b^1_j\otimes a^0_k \mapsto & -(y_j\cdot(x^i\be))(a)y^j(\al(c))x_i
	=-\be_{*a}(\al(c)) \\
	& =\rho(-\pmu(\pmu(\tilde c)\cdot \tilde a), -\pmu(\tilde c)\cdot \tilde a))\;,\\
b^0_j\otimes[b^1_j,a^1_k]\otimes a^0_k \mapsto & x^i(\be(a))y^j(\al(c))[x_i,y_j]
	=[\be(a),\al(c)] \\
	& =\rho(-[\pmu\tilde a,\pmu\tilde c],\pmu(\tilde c)\cdot \tilde a)\;,\\
b^0_j\otimes a^1_k\otimes b^1_j\cdot a^0_k \mapsto & x^i(\be(a))(x_i\cdot y^j\al)(c) y_j
	=\al_{*c}(\be(a)) \\
	& =\rho(\pmu(\pmu(\tilde a)\cdot \tilde c)+\pmu[\tilde a,\tilde c]),0)\;,\\
a^1_i\cdot b^0_l\otimes a^0_i\otimes b^1_l \mapsto & (x_i\cdot y^j\be)(a)x^i(\al(b)) y_j
	=\be_{*a}(\al(b)) \\
	& =\rho(\pmu(\pmu(\tilde b)\cdot \tilde a),\pmu(\tilde b)\cdot \tilde a))\;,\\
b^0_j\otimes b^1_j\cdot b^0_l\otimes b^1_l \mapsto & x^i(\be(a))(x_i\cdot y^j\be)(b) y_j
	=\be_{*b}(\be(a)) \\
	& =\rho(-\pmu(\pmu(\tilde a)\cdot \tilde b)-\pmu[\tilde a,\tilde b],
		-(\pmu(\tilde a)\cdot \tilde b)-[\tilde a,\tilde b])\;,\\
b^0_j\otimes b^0_l\otimes [b^1_j,b^1_l] \mapsto & x^i(\be(a)) y^j(\be(b)) [x_i,y_j]
	=[\be(a),\be(b)] \\
	& =\rho([\pmu(\tilde a),\pmu(\tilde b)],(\pmu(\tilde a)\cdot\tilde b)
		-(\pmu(\tilde b)\cdot\tilde a)+[\tilde a,\tilde b])\;,\\
\end{align*}
and we get \eqref{eq:cybe2} and \eqref{eq:cybe3} as a result of \eqref{eq:dba}.
\end{proof}

\section{Proof of Theorem \ref{tth}}

\begin{proof}
Let $r$ be a geometric solution of the CYBE. Theorem \ref{tle}
attaches to $r$ a classical bijective cocycle $7$-tuple
$(\lieg,\liea,X,\rhoga,\rhogax,\pi,\Psi_{\liea})$.  We will now
``exponentiate'' this classical 7-tuple to produce a (formal) quantum
bijective cocycle 7-tuple.

Recall that the category of (complex) Lie algebras is equivalent to the category
of (complex) formal groups, via the exponentiation functor $\lieh\to e^\lieh$.
The formal group $e^{\lieh}$ is a scheme which can be evaluated on
pro-Artinian local complex algebras; we will use, however, only
$e^{\lieh}(\CC[[\hbar]])=\{e^{\sum_{n>0}a_n\hbar^n}\ |\ a_n\in\lieh\ \forall n\}$,
and we drop $(\CC[[\hbar]])$ from the notation.
Thus we may consider the formal groups $G=e^{\lieg}$, $A=e^{\liea}$.

To exponentiate $\rhoga$,
we notice that $\rhoga:\lieg\to\Der(\liea)$ is a homomorphism.
Hence we have a homomorphism
$e^{\rhoga}:e^\lieg\to e^{\Der(\liea)}\subseteq\Aut(e^\liea)$.
For $\rhogax$, we have 
a homomorphism $\rhogax:\liea\rtimes\lieg\to\vect(X)$, 
hence we have a homomorphism
$e^{\rhogax}:e^{\liea\rtimes\lieg}\to e^{\vect(X)}=\Aut(X)$, where
$\Aut(X)$ stands for the group of formal diffeomorphisms $X\to X$.

We have the following short exact sequence of Lie algebras:
$0\to\liea\to\liea\rtimes\lieg\to\lieg\to 0$. By exponentiation, it maps to
$1\to e^\liea\to e^{\liea\rtimes\lieg}\to e^\lieg\to 1$. Since
the former sequence splits, the latter also does. This gives an isomorphism
$e^{\liea\rtimes\lieg}\simeq e^\liea\rtimes e^\lieg$. Using this isomorphism
we can consider $e^{\rhogax}$ to be a map $A\rtimes G\to\Aut(X)$. 

Consider now the map
$\pi:\lieg\to\liea$. This is a $1$-cocycle and hence yields a Lie algebra map
$\bar\pi:\lieg\to\liea\rtimes\lieg$, $x\mapsto (\pi(x),x)$. We now exponentiate
it and we get $e^{\bar\pi}:e^\lieg\to e^{\liea\rtimes\lieg}$ and via the previous
isomorphism we get a map $e^\lieg\to e^\liea\rtimes e^\lieg$. This is
a morphism of groups, and it must be of the form $e^g\mapsto (\widetilde\pi (e^g),e^g)$,
for some $\widetilde\pi:e^\lieg\to e^\liea$, which, a fortiori, is a $1$-cocycle.
Last, take $\widetilde\Psi_{\liea}=e^{\hbar\Psi_\liea}$, i.e.,
$\widetilde\Psi_{\liea}(x)=e^{\hbar\Psi_\liea(x)}$.

We have to prove that this is a bijective cocycle $7$-tuple. It is easy to see
that $\widetilde\Psi_\liea$ is $A\rtimes G$-equivariant; the rest of the conditions
are clear.

Now define $R$ by formula \eqref{eq:sos}.
From Lemma \ref{lm:sfst} we know that $R$ is a solution of the 
quantum Yang-Baxter equation,
and it is evident from the construction that it is geometric.
To see that $r$ is the classical limit of $R$, we compute the first order
approximation of \eqref{eq:sos}.
We can see that $x\check *y\equiv x-\hbar(y\aaad x)$ and
$x\check\circ y\equiv y-\hbar(y\acai x)$ modulo $\hbar^2$.
Therefore,
\begin{align*}
(R(f\otimes g))(x,y)&=(f\otimes g)R^{-1}(x,y)
	\equiv(f\otimes g)(x+\hbar(y\aaad x), y+\hbar(y\acai x))\ (\mbox{mod}\;\hbar^2) \\
	&=f(x)g(y)+\hbar\left(f(x)b^0_i(x)(b^1_i\cdot g)(y)+(a^1_j\cdot f)(x)a^0_j(y)g(y)\right),
\end{align*}
whence $R=1+\hbar r+\Oo(\hbar^2)$.
\end{proof}

\begin{remar}[Special cases]
Let us point out two special cases of this construction. The first one is
when $r$ is unitary ($r^{21}+r=0$). In this case 
our construction coincides with the one in \cite{es}; in particular,
we get $R$ unitary ($R^{21}R=1$).
On the CBCST side, the property of being unitary is equivalent to
$\liea$ being abelian, and it can be shown that when this happens
$\widetilde\pi(e^g)=e^{\frac{e^{g}-1}{g}\pi(g)}$ (this formula appears
in \cite{es}).
The other special case is
$r$ being a \emph{classical rack} ($r=b^0_j\otimes b^1_j$). In this case we obtain
a geometric rack, $R(x,y)=(x,x\acad y)$ for a suitable $\acad$ (see for instance
\cite{fr} for the definition of \emph{rack}).
\end{remar}

\section{Example}
In this section we apply the previous procedure to a $3$-dimensional example.
Let $X=\CC^3=\{(x_1,x_2,x_3)\}$, $r=r^1+r^2$, where
\;$r_1=\sum x_i\otimes A_i$,\; $r_2=\sum B_i\otimes x_i$.
Here $x_i$ stands for the canonical coordinate function and $A_i$, $B_i$ are the vector fields
defined by:
\begin{align*}
A_1(x,y,z)=(-\frac x2+y)\partial_x-\frac z2\partial_z, \qquad&
	B_1(x,y,z)=(-\frac x2+y)\partial_x+(y-\frac z2)\partial_z, \\
A_2(x,y,z)=-x\partial_x+y\partial_y, \qquad&
	B_2(x,y,z)=-x\partial_x+y\partial_y-x\partial_z, \\
A_3(x,y,z)=x\partial_x+z\partial_z, \qquad&
	B_3(x,y,z)=x\partial_x+z\partial_z.
\end{align*}

Following the definitions in \S\ref{ss2}, we see that $\liea$ can be identified
with the Heisenberg algebra
$$\liea=\spn\{X,Y,C\},\qquad \text{$C$ is in the center and $[X,Y]=C$},$$
and $\lieg$ can be identified with the upper-triangular matrices in $\mathfrak{gl}_2(\CC)$,
$$\lieg=\{\mdpd pq0r\}.$$
The map $\Psi:X\to\liea$ is just the ``identity"
$\Psi(x_1,x_2,x_3)=x_1X+x_2Y+x_3C$ and thus we will just identify $X$ with $\liea$.
The action $\rhoga$ is given by left multiplication by the matrix
$\rhoga(\mdpd pq0r)=\mtpt pq00r000{p+r}$ in the basis $\{X,Y,C\}$.
The $1$-cocycle is given by $\pi(\mdpd pq0r)=qX+rY+(p+\frac q2+r)C$.

In particular, the adjoint action of $\liea\rtimes\lieg$ on itself is given by
$$\adj(aX+bY+cC,\mdpd pq0r) = \left(\begin{smallmatrix}
p&q&0&-a&-b&0\\0&r&0&0&0&-b\\-b&a&p+r&-c&0&-c\\0&0&0&0&0&0\\0&0&0&-q&p-r&q\\0&0&0&0&0&0
\end{smallmatrix}\right)$$
Considering the adjoint representation, which is faithful, we can compute $\widetilde\pi$ as
$$\widetilde\pi(e^{\hbar\mdpd pq0r})
	=(e^{\adj})^{-1}\left(e^{\hbar\adj(\pi\mdpd pq0r)+\hbar\adj\mdpd pq0r}e^{-\hbar\adj\mdpd pq0r}\right).$$
From this we get that $\widetilde\pi^{-1}(e^{\hbar(aX+bY+cC)})=e^{\mdpd pq0r}$, where
\begin{align*}
p &=\ln(\frac{1+\hbar c-\frac 12\hbar a}{1+\hbar b}), \\
q &=\hbar a \left(1+\hbar b\right) \ln\left(\frac {1+\hbar c-\frac 12\hbar a}{\left(1+\hbar b\right)^2}\right)
	\left(1+\hbar c-\frac 12\hbar a-\left(1+\hbar b\right)^2\right)^{-1}, \\
r &=\ln(1+\hbar b).
\end{align*}
This computation, as well as most of the remaining ones, were carried with the help
of MAPLE.
Finally, we can compute in a straightforward way $R(x,y)=(x\check *y, x\check\circ y)$
as defined in \eqref{eq:sos}, and we get
\begin{align*}
(x_1,x_2,x_3) \check * (y_1,y_2,y_3)
	=\Big(&\frac{1-\hbar y_3+\hbar \frac{y_1}2}{1-\hbar y_2}x_1-\hbar y_1x_2,\;
		(1-\hbar y_2)x_2,\;
		(1-\hbar y_3+\hbar \frac{y_1}2) x_3\Big); \\
(x_1,x_2,x_3) \check \circ (y_1,y_2,y_3) =\Big(
 & \frac{y_1(1-\hbar x_2)
	+\hbar y_2(x_1 - y_1 + \hbar y_1x_2 - \hbar x_1y_3 + \frac 12\hbar x_1y_1)}{\text{den}}, \\
 & \frac{y_2}{(1-\hbar x_2+\hbar^2x_2y_2)}, \\
 & \frac{(1-\hbar y_2)(y_3 - \hbar y_1x_2)
	+\hbar y_2x_1(1 - \hbar y_3 + \frac 12\hbar y_1)}{\text{den}}\Big)
\end{align*}
where
$$\text{den}=
(1 - \hbar y_2)(1 - \frac 12\hbar^2 y_1x_2)
	+\hbar(1 - \hbar y_3 + \frac 12\hbar y_1)(-x_3 + \frac 12x_1 + \hbar x_3y_2).$$

It is possible to ``unitarize" this example by replacing in $\liea$ the bracket
$[X,Y]=C$ by $[X,Y]=\ve C$. The r-matrix $r^\ve$ has a similar expresion to that
of $r$, but changing $A_i$, $B_i$ by
\begin{align*}
A_1^\ve(x,y,z)=(-\ve\frac x2+y)\partial_x-\ve\frac z2\partial_z, \qquad&
	B_1^\ve(x,y,z)=(-\ve\frac x2+y)\partial_x+\ve(y-\frac z2)\partial_z, \\
A_2^\ve(x,y,z)=-x\partial_x+y\partial_y, \qquad&
	B_2^\ve(x,y,z)=-x\partial_x+y\partial_y-\ve x\partial_z, \\
A_3^\ve(x,y,z)=x\partial_x+z\partial_z, \qquad&
	B_3^\ve(x,y,z)=x\partial_x+z\partial_z.
\end{align*}
The $1$-cocycle will become
$\pi^\ve(\mdpd pq0r)=qX+rY+(p+\frac{\ve q}2+r)C$.
There are similar expresions for the R-matrix $R^\ve$. Now, if we let $\ve\to 0$, we will get
a unitary r-matrix (since $\liea$ will become abelian) and hence a unitary R-matrix.
In the limit, we get
$\lim_{\ve\to 0}R^\ve(x,y)=(x\check *^0 y, x\check\circ^0 y)$, given by
\begin{align*}
(x_1,x_2,x_3) \check *^0 (y_1,y_2,y_3)
=\Big(&
\frac{1-\hbar y_3}{1-\hbar y_2}x_1-\hbar y_1x_2, (1-\hbar y_2)x_2, (1-\hbar y_3)x_3\Big) \\
(x_1,x_2,x_3) \check \circ^0 (y_1,y_2,y_3)
=\Big(&
y_1\frac{1-\hbar x_2}{1-\hbar x_3+\hbar ^2x_3y_3}
	+\hbar y_2\frac{x_1(1-\hbar y_3)}{(1-\hbar x_3+\hbar ^2x_3y_3)(1-\hbar y_2)}, \\
&\frac{y_2}{1-\hbar x_2+\hbar^2x_2y_2}, \\
&\frac{y_3}{(1-\hbar x_3+\hbar ^2x_3y_3)}\Big)
\end{align*}

\end{document}